\documentclass[12pt,english]{article}
\usepackage[latin9]{inputenc}
\usepackage{geometry}
\usepackage{amsmath}
\usepackage{amssymb}

\makeatletter
\date{}

\makeatother

\usepackage{babel}
\begin{document}

\title{Coordinate deletion of zeroes}

\author{Eero R{\"a}ty\thanks{Centre for Mathematical Sciences, Wilberforce Road, Cambridge CB3
0WB, UK, epjr2@cam.ac.uk}}
\maketitle
\begin{abstract}
For a family $A\subseteq\left\{ 0,\dots,k\right\} ^{n}$, define the
$\delta$-shadow of $A$ to be the set obtained from $A$ by removing
from any of its vectors one coordinate that equals zero. Given the
size of $A$, how should we choose $A$ to minimise its $\delta$-shadow?
Our aim in this paper is to show that, for any $r$, the family of
all sequences with at most $r$ zeros has minimal $\delta$-shadow.
We actually give the exact best $A$ for every size. 
\end{abstract}

\section{Introduction}

The classical Kruskal-Katona theorem is concerned with the lower shadow
of set systems. For $A\subseteq\left\{ 0,1\right\} ^{n}$, define
the \textit{lower shadow} of $A$ to be the set of sequences obtained
from any of its vectors by flipping one of its 1-entries to 0. The
\textit{rank} of a sequence $x\in\left\{ 0,\dots,n\right\} ^{k}$
is defined to be $\left|x\right|=\sum_{i=1}^{k}x_{i}$. Note that
the lower shadow operator decreases the rank of a sequence by 1. For
given $r$, it is natural to ask how to choose a family $A\subseteq\left\{ 0,1\right\} ^{n}$
of given size containing only vectors with rank $r$, which minimises
the lower shadow. This question was answered by Kruskal \cite{key-5}
and Katona \cite{key-4}.

Define the \textit{colexicographic order} on $\left\{ x\in\left\{ 0,1\right\} ^{n}\,:\,\left|x\right|=r\right\} $
by $x\leq_{colex}y$ if $\max\left(X\Delta Y\right)\in Y$. Here $X=\left\{ i\,:\,x_{i}=1\right\} $
and $Y=\left\{ i\,:\,y_{i}=1\right\} $ . The Kruskal-Katona theorem
states that for a set $A\subseteq\left\{ 0,1\right\} ^{n}$ containing
only sequences of rank $r$, the lower shadow is minimised when $A$
is chosen to be an initial segment of colexicographic order.

Instead of changing the coordinates, it is also natural to define
an operator which acts by deleting coordinates. For $A\subseteq\left\{ 0,\dots,k\right\} ^{n}$
define the \textit{coordinate deletion shadow} $\Delta A$ to be the
set of those sequences obtained from any of its vectors by deleting
one coordinate. For example $\Delta\left(\left\{ 000,001,002,121\right\} \right)=\left\{ 00,01,02,12,11,21\right\} $. 

Again it is natural to ask that which sets minimises the coordinate
deletion shadow. Define the \textit{simplicial order} $\leq_{sim}$
on $\left\{ 0,1\right\} ^{n}$ by 
\[
x\leq_{sim}y\text{ if }\left|x\right|<\left|y\right|\text{ or }\left|x\right|=\left|y\right|\text{ and }\min(X\Delta Y)\in X.
\]

It was proved by Danh and Daykin that for subsets of $\left\{ 0,1\right\} ^{n}$,
$\Delta A$ is minimised by an initial segment of the simplicial order
\cite{key-3}. They also conjectured a certain order as best in $\left\{ 0,1,\dots\right\} ^{n}$,
but Leck \cite{key-6} showed that this turned out to be false and
in fact there is no order in general whose all initial segments have
minimal coordinate deletion shadow. 

Bollob{\'a}s and Leader \cite{key-1} pointed out that for $k\geq2$ the
sets $A_{t}=\left\{ 0,\dots,t-1\right\} ^{n}\subseteq\left\{ 0,\dots,k\right\} ^{n}$
are extremal for $\Delta$. Indeed, suppose that $B\subseteq\left\{ 0,\dots,k\right\} ^{n}$
is extremal with $\left|B\right|=t^{n}$. Define $B_{\left[n\right]\setminus\left\{ i\right\} }$
to be the projection of $B$ onto the hyperplane excluding the $i^{th}$
direction. Suppose that $\left|\Delta B\right|<t^{n-1}$. Since $B_{\left[n\right]\setminus\left\{ i\right\} }\subseteq\Delta B$
for all $i$, it follows that $\left|B_{\left[n\right]\setminus\left\{ i\right\} }\right|<t^{n-1}$
for all $i$. Thus the Loomis-Whitney inequality \cite{key-7} implies
that $\left|B\right|^{n-1}\leq\left(\prod_{i=1}^{n}B_{\left[n\right]\setminus\left\{ i\right\} }\right)<t^{n\left(n-1\right)}$,
which contradicts $\left|B\right|=t^{n}$. Hence $\left|\Delta B\right|\geq t^{n-1}$
and since $\Delta A_{t}=t^{n-1}$, it follows that each $A_{t}$ is
extremal. 

Bollob{\'a}s and Leader also made the following conjecture that certain
other type of sets are also extremal.\\

\textbf{Conjecture 1 }\cite{key-1}\textbf{.} For each $t\leq k$
and $r\leq k$, let $B_{r,t}\subseteq\left\{ 0,\dots,k\right\} ^{n}$
be the subset containing all sequences with at most $r$ zeros, and
with all coordinates in $\left\{ 0,\dots,t\right\} $. Then the sets
$B_{r,t}$ have extremal $\Delta$-shadow in $\left\{ 0,\dots,k\right\} ^{n}$.\\
\\
Even the case $t=k$ in the conjecture in unknown. 

There is, however, a notion that comes 'between' the lower shadow
and the coordinate deletion shadow. The usual lower shadow operator
decreases the rank by 1 and preserves the dimension $n$, while the
coordinate deletion shadow decreases the dimension by 1 but there
is no control on how it changes the rank. So it is natural to consider
the following operator which preserves the rank, but reduces the dimension
by one.

Define the $\delta$\textit{-shadow} of $A\subseteq\left\{ 0,\dots,k\right\} ^{n}$
to be the set of sequences in $\left\{ 0,\dots,k\right\} ^{n-1}$
obtained by removing one coordinate that equals 0 from any of the
vectors in $A$. Denote this set by $\delta A$. Thus for example
$\delta\left(\left\{ 00011,00101\right\} \right)=\left\{ 0011,0101\right\} $
and $\delta\left(\left\{ 112,113,123\right\} \right)=\emptyset$.

How can we find sets $A$ with minimal $\delta$-shadow? If $\left|A\right|\leq k^{n}$
then the question is trivial, as one can take any subset of $\left\{ 1,\dots,k\right\} ^{n}$
of given size. In general, it is natural to choose $A$ to contain
sequences with as few zeros as possible. Furthermore, it is natural
to guess that for each $0\leq i\leq n$, the sets containing all sequences
with at most $i$ zeros have minimal $\delta$-shadow. 

Our main result in this paper is to find an order on $\left\{ 0,\dots,k\right\} ^{n}$
whose initial segments have minimal $\delta$-shadow. In particular,
it follows that the sets containing all sequences with at most $i$
zeros have minimal $\delta$-shadow. 

In order to state the main result, we need a few definitions. For
a sequence $x\in\left\{ 0,\dots,k\right\} ^{n}$, let $R\left(x\right)=\left\{ i\,:\,x_{i}=0\right\} $
and let $w\left(x\right)=\left|R\left(x\right)\right|$. Let $L_{r}\left(n\right)=\left\{ x\in\left\{ 0,\dots,k\right\} ^{n}\,:\,w\left(x\right)=r\right\} $.
Note that $\delta$ maps sequences in $L_{r}\left(n\right)$ to sequences
in $L_{r-1}\left(n-1\right)$.

For $x\in\left\{ 0,\dots,k\right\} ^{n}$, define its \textit{reduced
sequence} to be the sequence obtained by removing all coordinates
from $x$ that equal 0. Denote the reduced sequence of $x$ by $re(x)$.
Note that for any sequence $s$ and for any $t\in\delta s$ we have
$re(s)=re(t)$, as removing a coordinate which equals 0 does not change
the reduced sequence. Hence each $L_{r}\left(n\right)$ splits into
disjoint components based on the reduced sequences. 

We will start by proving that inside a component one should choose
sequences $x$ for which the sets $R\left(x\right)$ form an initial
segment of colex. This is a straightforward consequence of the work
of Danh and Daykin in \cite{key-3}. 

Since we know that $\left\{ 0,\dots,k\right\} ^{n}$ splits into components
based on the reduced sequences, and we know that initial segments
of the colexicographic order minimises the $\delta$-shadow inside
each component, we are left with the question on how to split the
sequences into different components in order to minimise the shadow. 

We go on to prove that in order to minimise the shadow of a subset
in $\left\{ 0,\dots,k\right\} ^{n}$, one should first choose sequences
in components in $L_{r}\left(n\right)$ rather than in $L_{s}\left(n\right)$
for all $r<s$, and inside $L_{r}\left(n\right)$ one should choose
all sequences from a component before taking any sequences from another
component. As a consequence we obtain an order whose initial segments
minimises the $\delta$-shadow. 

For $r\in\left\{ 0,1,\dots,k\right\} $ define $R_{r}\left(x\right)=\left\{ i\,:\,x_{i}=r\right\} $
and $w_{r}\left(x\right)=\left|R_{r}\left(x\right)\right|$. Note
that $R=R_{0}$ and $w=w_{0}$. For all $r$ we define an order $\leq_{c}$
on $\left\{ 1,\dots,k\right\} ^{r}$ as follows. For distinct $x$,
$y\in\left\{ 1,\dots,k\right\} ^{r}$ let $i$ be minimal such that
$R_{i}\left(x\right)\neq R_{i}\left(y\right)$. We say that $x\leq_{c}y$
if and only if $\max\left(R_{i}\left(x\right)\Delta R_{i}\left(y\right)\right)\in R_{i}\left(y\right)$. 

Define an order $\leq$ on $\left\{ 0,\dots,k\right\} ^{n}$ as follows.
For distinct $x$, $y\in\left\{ 0,\dots,k\right\} ^{n}$ we set $x\leq y$
if 
\begin{enumerate}
\item $w_{0}\left(x\right)<w_{0}\left(y\right)$
\item $w_{0}\left(x\right)=w_{0}\left(y\right)$, $re\left(x\right)\neq re\left(y\right)$
and $re\left(x\right)\leq_{c}re\left(y\right)$
\item $w_{0}\left(x\right)=w_{0}\left(y\right)$, $re\left(x\right)=re\left(y\right)$
and $R_{0}\left(x\right)\leq_{colex}R_{0}\left(y\right)$
\end{enumerate}
Now we are ready to state our main theorem. \\

\textbf{Theorem 1. }Let $A\subseteq\left\{ 0,\dots,k\right\} ^{n}$
and let $B$ be an initial segment of $\leq$ with $\left|B\right|=\left|A\right|$.
Then $\left|\delta A\right|\geq\left|\delta B\right|$.\\

In particular, it follows that the sets of the form $L_{\leq r}\left(n\right)=\bigcup_{i=0}^{r}L_{i}\left(n\right)$
are extremal. Note that for fixed $r$, every component of $L_{r}\left(n\right)$
behaves in the same way. Hence for any fixed $r$, one could replace
the $\leq_{c}$ order by any other order on $\left\{ 1,\dots,k\right\} ^{r}$
in the definition of the $\leq$-order. 

The plan of the paper is as follows. In Section 2 we prove that inside
a component the sets of sequences whose associated sets $R\left(x\right)$
form an initial segment of colex have minimal $\delta$-shadow. In
Section 3 we prove Theorem 1. In Section 4 we generalise the $\delta$-shadow
to allow deleting coordinates in some set $\left\{ 0,\dots,r\right\} $
instead of just deleting only coordinate which equals 0. In this case
we will show that sets $\left\{ x\,:\,\sum_{i=0}^{r}w_{i}\left(x\right)\leq s\right\} $,
which are analogous to the sets $L_{\leq s}\left(n\right)$, are extremal
for all $0\leq s\leq n$. In this general case we do not know what
happens for sets of other sizes. 

We write use the standard notation $\left[n\right]=\left\{ 1,\dots,n\right\} $
and $\left[n\right]^{(r)}=\left\{ A\subseteq\left[n\right]\,:\,\left|A\right|=r\right\} $.
We write $L_{r}$ instead of $L_{r}\left(n\right)$ if the dependence
on $n$ is clear. When $k=1$ we may also write $\left\{ 0,1\right\} _{r}^{n}$
instead of $L_{r}\left(n\right)$. This notation will be used to highlight
that we are working with $\left\{ 0,1\right\} $-sequences. 

\section{Deletion on $\left\{ 0,1\right\} $-sequences}

In \cite{key-3} Danh and Daykin proved the following result for the
coordinate deletion shadow $\Delta$ on $\{0,1\}^{n}$.\\

\textbf{Theorem 2 (Danh, Daykin). }Let $A\subseteq\left\{ 0,1\right\} ^{n}$
and let $B$ be an initial segment of simplicial order with $\left|B\right|=\left|A\right|$.
Then $\left|\Delta A\right|\geq\left|\Delta B\right|$. $\square$\\

There is a natural correspondence between the sequences $\left\{ 0,1\right\} ^{n}$
and the power-set ${\cal P}\left(\left\{ 1,\dots,n\right\} \right)$.
For our purposes it will be convenient to choose this correspondence
to be given by mapping a sequence $\left(x_{i}\right)$ to the set
$R_{0}\left(x\right)=\left\{ i\,:\,x_{i}=0\right\} $. 

In this way we can identify set $A\subseteq L_{r}\left(n\right)$
with a set system ${\cal A}\subseteq\left[n\right]^{\left(r\right)}$
by taking ${\cal A}$ to be the images of the elements of $A$ under
this bijection. This enables us to translate questions on $\delta$
to questions related to properties of set systems ${\cal A}\subseteq\left[n\right]^{\left(r\right)}$
instead. We start by proving that the subsets $A$ of $L_{r}$ with
minimal shadow are the ones whose corresponding set ${\cal A}$ is
an initial segment of colex. \\

\textbf{Lemma 3. }Let $A\subseteq\left\{ 0,1\right\} _{r}^{n}$, and
let $B\subseteq\left\{ 0,1\right\} _{r}^{n}$ be a set of same size
for which ${\cal B}$ is an initial segment of colex. Then $\left|\delta A\right|\geq\left|\delta B\right|$.\\

\textbf{Proof. }Define $C_{1}=A\cup L_{>r}\left(n\right)$ and $C_{2}=B\cup L_{>r}\left(n\right)$,
where $L_{>r}\left(n\right)=\bigcup_{i=r+1}^{n}\left\{ 0,1\right\} _{i}^{n}$
. Now $C_{2}$ is isomorphic to initial segment of simplicial order,
and the isomorphism is the map which reverses the sequences. Since
this map preserves the size of $\Delta$-shadow, Theorem 2 implies
that $\Delta C_{2}$ is minimal and hence 
\begin{equation}
\left|\Delta C_{2}\right|\leq\left|\Delta C_{1}\right|.
\end{equation}
Note that $\Delta C_{1}=L_{>r}(n-1)\cup\delta A$ and $\Delta C_{2}=L_{>r}(n-1)\cup\delta B$.
Indeed, $L_{>r}(n-1)$ is certainly subset of both of these sets,
and the only contribution to elements not in $L_{>r}(n-1)$ comes
from removing 0 from a sequence which contains exactly $n-r$ 1's.
Hence 
\begin{equation}
\left|\Delta C_{1}\right|=\left|L_{>r}\left(n-1\right)\right|+\left|\delta A\right|
\end{equation}
and 
\begin{equation}
\left|\Delta C_{2}\right|=\left|L_{>r}\left(n-1\right)\right|+\left|\delta B\right|.
\end{equation}
Combining (1), (2) and (3) yields that $\left|\delta A\right|\ge\left|\delta B\right|$.
$\square$\\

Lemma  3 implies that inside $\left\{ 0,1\right\} _{r}^{n}$ the colexicographic
order minimises the size of the shadow. Before moving on to general
$k$ from $k=1$, we find a way to relate the size of $\delta A$
to the associated family ${\cal A}$. For convenience, from now on
we say that $A\subseteq\left\{ 0,1\right\} _{r}^{n}$ is an initial
segment of colex if the associated set system ${\cal A}$ is an  initial
segment of colex. For ${\cal A}\subseteq{\cal P}\left(\left\{ 1,\dots,n\right\} \right)$
define ${\cal A}_{1}=\left\{ B\in{\cal A}\,:\,1\in B\right\} $.\\

\textbf{Lemma 4. }Let $A\subseteq\left\{ 0,1\right\} _{r}^{n}$ be
an initial segment of colex associated to ${\cal A}$. Then $\left|\delta A\right|=\left|{\cal A}_{1}\right|$.\\

\textbf{Proof. }The proof is by induction on $\left|A\right|$, the
case $\left|A\right|=1$ is clear. Let $B$ be an initial segment
with $\left|B\right|=\left|A\right|+1$, say $B=A\cup\left\{ x\right\} $
with $x=x_{1}\dots x_{n}$. First we will prove that $x_{2}\dots x_{n}$
is the only element which could be in $\delta B\setminus\delta A$. 

Indeed, suppose $t\in\delta B\setminus\delta A$ and that it is obtained
by removing the $k^{th}$ coordinate of $x$. Hence $t=x_{1}\dots x_{k-1}x_{k+1}\dots x_{n}$
and $x_{k}=0$. Let $i=\min\left\{ j\,:\,x_{j}=1\right\} $ and set
$y=0t=0x_{1}\dots x_{k-1}x_{k+1}\dots x_{n}$. If $i\leq k$, then
$y_{j}=x_{j}$ for all $j\leq i-1$ but $y_{k}=x_{k-1}=0\neq1=x_{k}$
so $y<_{colex}x$. But then $t\in\delta y\subseteq\delta A$ which
contradicts $t\in\delta B\setminus\delta A$. 

Hence we must have $i>k$. But in this case $x_{1}=\dots x_{k}=0$
and therefore $t=0\dots0x_{k+1}\dots x_{n}=x_{2}\dots x_{n}$. Hence
$\delta B\setminus\delta A$ is either empty or contains only $x_{2}\dots x_{n}$.

Note that $0x_{2}\dots x_{n}$ is the least element in colex which
has $x_{2}\dots x_{n}$ contained in its $\delta$-shadow. Hence $0x_{2}\dots x_{n}\in\delta B\setminus\delta A$
if and only if $x=0x_{2}\dots x_{n}$. Thus 
\[
\left|\delta B\right|=\begin{cases}
\begin{array}{c}
\left|\delta A\right|+1\\
\left|\delta A\right|
\end{array} & \begin{array}{c}
\text{if }x_{1}=0\\
\text{if }x_{1}=1
\end{array}\end{cases}.
\]
Also ${\cal B}={\cal A}\cup R_{0}(x)$, and the set $R_{0}(x)$ contains
$1$ if and only if $x_{1}=0$. Thus 
\[
\left|{\cal B}_{1}\right|=\begin{cases}
\begin{array}{c}
\left|{\cal A}_{1}\right|+1\\
\left|{\cal A}_{1}\right|
\end{array} & \begin{array}{c}
\text{if }x_{1}=0\\
\text{if }x_{1}=1
\end{array}\end{cases}
\]
and hence $\left|\delta B\right|=\left|{\cal {\cal B}}_{1}\right|$
by induction. $\square$

\section{The main theorem}

Let $H$ be the bipartite graph with vertex set $\left\{ 0,\dots,k\right\} ^{n}\cup\left\{ 0,\dots,k\right\} ^{n-1}$
and whose edges are precisely those pairs $s,\,t$ with $s\in\left\{ 0,\dots,k\right\} ^{n}$
and $t\in\delta s$. Then for $A\subseteq\left\{ 0,\dots,k\right\} ^{n}$,
$\delta A$ is just the neighbourhood of $A$ in the graph $H$. Note
that both classes can be partitioned as $\left\{ 0,\dots,k\right\} ^{n}=\bigcup_{i=0}^{n}L_{i}\left(n\right)$
and $\left\{ 0,\dots,k\right\} ^{n-1}=\bigcup_{i=0}^{n}L_{i}\left(n-1\right)$,
and by definition of $\delta$ it is clear that there are edges only
between $L_{i}\left(n\right)$ and $L_{i-1}\left(n-1\right)$, with
the convention $L_{-1}=\emptyset$.

Let $C$ be a connected component in $H$. Suppose $C$ is non-trivial,
i.e. $\left\{ 0,\dots,k\right\} ^{n}\cap C\subseteq L_{i}(n)$ for
some $i>0$. Recall that for all $x$ and for any $y\in\delta x$,
$x$ and $y$ have the same reduced sequences. But since $C$ is a
connected component, this means that every $x\in C$ has the same
reduced sequence. Conversely it is easy to check that for $i>0$ all
sequences $x\in L_{i}(n)\cup L_{i-1}(n-1)$ with the same reduced
sequence are in the same connected component. Thus we can deduce that
the connected components in $H$ are given as follows.\\

\textbf{Lemma 5. }For $s\in\bigcup_{i=0}^{r}\left\{ 1,\dots,k\right\} ^{i}$
define $C_{s}=\left\{ x\in\left\{ 0,\dots,k\right\} ^{n}\,:\,re(x)=s\right\} $
and $D_{s}=\left\{ x\in\left\{ 0,\dots,k\right\} ^{n-1}\,:\,re(x)=s\right\} $.
Then $C_{s}\cup D_{s}$ are the connected components of $H$. $\square$\\

Broadly speaking, we need to only understand how to minimise $\delta$
inside a connected component and to determine how to distribute the
sequences into different connected components in order to minimise
$\delta$. It turns out that inside connected component one should
choose sequences $x$ whose sets $R_{0}\left(x\right)$ forms initial
segment of colex. \\

\textbf{Lemma 6. }Let $C\subseteq L_{i}(n)\cup L_{i-1}(n-1)$ be a
connected component corresponding to a reduced word $x=x_{1}\dots x_{n-i}$.
Let $B\subseteq L_{i}\cap C$ and let $A\subseteq L_{i}\cap C$ be
a set of sequences chosen such that $\left|A\right|=\left|B\right|$
and $\left\{ R_{0}\left(x\right)\,:\,x\in A\right\} $ is an initial
segment of colex. Then $\left|\delta B\right|\geq\left|\delta A\right|$.\\

Proof Note that the behaviour of the connected component depends only
on $n-i$ and in particular not on the sequence $x_{1}\dots x_{n-i}$,
as the reduced sequence and the order of coordinates in the reduced
sequence is preserved under taking $\delta$-shadow. In particular,
all such connected components has the same size and they all behave
in the same way under taking $\delta$-shadow. Hence it suffices to
consider only the component with $x_{1}=\dots=x_{n-i}=1$. But this
component is just $\left\{ 0,1\right\} _{i}^{n}$ and hence the result
follows from Lemma 3. $\square$\\

Hence it remains to understand how to fill different connected components.
Our aim is to show that it is optimal to first choose all sequences
in a component before taking sequences from another component, and
also to prefer a component in $L_{i}(n)$ over a component in $L_{i+1}(n)$. 

From now on we call the sets $C_{s}$ connected components, i.e. by
a connected component we refer to the intersection of a connected
component with $\left\{ 0,\dots,k\right\} ^{n}$. 

For $s,\,t\in\bigcup_{i=0}^{r}\left\{ 1,\dots,k\right\} ^{i}$ define
the \textit{$s,t$-compression} operator as follows. For $A\subseteq\left\{ 0,\dots,k\right\} ^{n}$
its compression $B=C_{s,t}\left(A\right)$ is given by setting 
\begin{enumerate}
\item $B\cap C_{s}$ to be an initial segment of colex of length $\min\left(\left|A\cap\left(C_{s}\cup C_{t}\right)\right|,\left|C_{s}\right|\right)$
\item $B\cap C_{t}$ to be an initial segment of colex of length $\max\left(0,\left|A\cap\left(C_{s}\cup C_{t}\right)\right|-\left|C_{s}\right|\right)$
\item $B\setminus\left(C_{s}\cup C_{t}\right)=A\setminus\left(C_{s}\cup C_{t}\right)$
\end{enumerate}
It is clear that $\left|C_{s.t}\left(A\right)\right|=\left|A\right|$
for all $s$ and $t$. As usual we say that $A\subseteq\left\{ 0,\dots,k\right\} ^{n}$
is \textit{$s,t$-compressed} if $C_{s,t}\left(A\right)=A$. 

In order to prove Theorem 1, we will need the following two Lemmas.\\

\textbf{Lemma 7. }Let $A\subseteq\left\{ 0,\dots,k\right\} ^{n}$
be a set and let $s,\,t\in\left\{ 1,\dots,k\right\} ^{n-i}$ for some
$i$. Then $\left|\delta A\right|\geq\left|\delta C_{s,t}\left(A\right)\right|$.
\\

\textbf{Lemma 8.} Let $A\subseteq\left\{ 0,\dots,k\right\} ^{n}$
be a set and let $s\in\left\{ 1,\dots,k\right\} ^{n-i}$, $t\in\left\{ 1,\dots,k\right\} ^{n-i-1}$
for some $i$. Then $\left|\delta A\right|\geq\left|\delta C_{s,t}\left(A\right)\right|$.\\

In order to prove these Lemmas, we will relate them to the appropriate
questions on the subsets of $\left[n\right]^{(i)}$. We will now state
these results, but the proof is presented after the proofs of Lemma
7 and Lemma 8.

Define ${\cal B}\subseteq\left[n\right]^{(r)}$ to be a \textit{segment}
if there exists initial segments ${\cal I}$ and ${\cal J}$ of colex
such that ${\cal A}={\cal I}\setminus{\cal J}$. \\

\textbf{Lemma 9. }The following claims are true.\\

\textbf{Claim 1.} Let ${\cal A}\subseteq\left[n\right]^{(i)}$ be
a segment and ${\cal I}\subseteq\left[n\right]^{(i)}$ be an initial
segment of colex with $\left|{\cal A}\right|=\left|{\cal I}\right|$.
Then $\left|{\cal A}_{1}\right|\leq\left|{\cal I}_{1}\right|$\\

\textbf{Claim 2.} Let ${\cal I}\subseteq\left[n\right]^{(i)}$ and
${\cal J}\subseteq\left[n\right]^{(i+1)}$ be initial segments of
colex with $\left|{\cal I}\right|=\left|{\cal J}\right|$. Then $\left|{\cal I}_{1}\right|\leq\left|{\cal J}_{1}\right|$\\

\textbf{Claim 3. }Let ${\cal A}\subseteq\left[n\right]^{(r)}$ be
a segment and let ${\cal I}=\left[n\right]^{(r)}\setminus{\cal J}$,
where ${\cal J}$ is an initial segment of colex chosen such that
$\left|{\cal I}\right|=\left|{\cal A}\right|$. Then $\left|{\cal A}_{1}\right|\geq\left|{\cal I}_{1}\right|$.\\

\textbf{Claim 4.} Let ${\cal I}_{*}$ and ${\cal J}_{*}$ be initial
segments of colex chosen such that ${\cal I}=\left[n\right]^{(i)}\setminus{\cal I}_{*}$
and ${\cal J}=\left[n\right]^{(i+1)}\setminus{\cal J}_{*}$ satisfies
$\left|{\cal I}\right|=\left|{\cal J}\right|$. Then $\left|{\cal I}_{1}\right|\leq\left|{\cal J}_{1}\right|$.\\

\textbf{Proof of Lemma 7.}

Let $A\subseteq\left\{ 0,\dots,k\right\} ^{n}$ and $B=C_{s,t}\left(A\right)$.
Note that $B$ depends only on $\left|A\cap C_{s}\right|$ and $\left|A\cap C_{t}\right|$.
Lemma 6 implies that an initial segment of colex minimises the $\delta$-shadow
inside a connected component, so we may assume that $Q=A\cap C_{s}$
and $R=A\cap C_{t}$ are initial segments of colex. 

Let $S=B\cap C_{s}$ and $T=B\cap C_{t}$. Let ${\cal Q}$, ${\cal R}$,
${\cal S}$ and ${\cal T}$ be the associated families in $\left[n\right]^{(i)}$.
Since $B\setminus\left(C_{s}\cup C_{t}\right)=A\setminus\left(C_{s}\cup C_{t}\right)$,
it follows that $\left|\delta B\right|\geq\left|\delta A\right|$
is equivalent to $\left|\delta Q\right|+\left|\delta R\right|\geq\left|\delta S\right|+\left|\delta T\right|$.
By applying Lemma 4, this can be rewritten as $\left|{\cal Q}_{1}\right|+\left|{\cal R}_{1}\right|\geq\left|{\cal S}_{1}\right|+\left|{\cal T}_{1}\right|$.\\

\textbf{Case 1. $\left|Q\right|+\left|R\right|\leq\left|C_{s}\right|$}\\

By definition of $B$, it follows that $T=\emptyset$ and $\left|S\right|=\left|Q\right|+\left|R\right|$.
Let ${\cal I}={\cal S}\setminus{\cal Q}$. Since ${\cal S}$ and ${\cal Q}$
are initial segments of colex, it follows that ${\cal I}$ is a segment
of length $\left|{\cal R}\right|$. Thus $\left|{\cal I}_{1}\right|\leq\left|{\cal R}_{1}\right|$
by Claim 1 and hence 
\[
\left|{\cal S}_{1}\right|+\left|{\cal T}_{1}\right|=\left|{\cal Q}_{1}\right|+\left|{\cal I}_{1}\right|\leq\left|{\cal Q}_{1}\right|+\left|{\cal R}_{1}\right|
\]
as required. \\

\textbf{Case 2. $\left|Q\right|+\left|R\right|>\left|C_{s}\right|$}\\

In this case $S=C_{s}$ and hence $|T|<|R|$. Thus we can write ${\cal I}={\cal R}\setminus{\cal T}$,
which is a segment as ${\cal R}$ and ${\cal T}$ are initial segments
of colex. Also set ${\cal J}={\cal S}\setminus{\cal Q}=\left[n\right]^{(i)}\setminus{\cal Q}$,
which is a segment as well. Since $\left|{\cal S}\right|+\left|{\cal T}\right|=\left|{\cal R}\right|+\left|{\cal Q}\right|$
it follows that $\left|{\cal I}\right|=\left|{\cal J}\right|$. Thus
Claim 3 implies that $\left|{\cal J}_{1}\right|\leq\left|{\cal I}_{1}\right|$.

Combining this together with the definitions of ${\cal I}$ and ${\cal J}$
implies that 
\[
\left|{\cal Q}_{1}\right|+\left|{\cal R}_{1}\right|=\left|{\cal Q}_{1}\right|+\left|{\cal I}_{1}\right|+\left|{\cal T}_{1}\right|\leq\left|{\cal Q}_{1}\right|+\left|{\cal J}_{1}\right|+\left|{\cal T}_{1}\right|=\left|{\cal S}_{1}\right|+\left|{\cal T}_{1}\right|
\]
as required, which completes the proof of Lemma 7. $\square$ \\

\textbf{Proof of Lemma 8.}

Let $A\subseteq\left\{ 0,\dots,k\right\} ^{n}$ and $B=C_{s,t}\left(A\right)$.
By Lemma 5 we may assume that both $A\cap C_{s}$ and $A\cap C_{t}$
are initial segments of colex. As in the proof of Lemma 7, set $Q=A\cap C_{s}$
, $R=A\cap C_{t}$, $S=B\cap C_{s}$ and $T=B\cap C_{t}$. Let ${\cal Q}$
and ${\cal S}$ be the associated set systems in $\left[n\right]^{(i)}$,
and ${\cal R}$ and ${\cal T}$ be the associated set systems in $\left[n\right]^{(i+1)}$.
By Lemma 4 it suffices to prove that $\left|{\cal Q}_{1}\right|+\left|{\cal R}_{1}\right|\geq\left|{\cal S}_{1}\right|+\left|{\cal T}_{1}\right|$.\\

\textbf{Case 1. $\left|{\cal Q}\right|+\left|{\cal R}\right|\leq\left|C_{s}\right|$}

By definition of $B$, it follows that ${\cal S}$ is an initial segment
of colex of length $\left|{\cal Q}\right|+\left|{\cal R}\right|$
in $\left[n\right]^{(i)}$, and ${\cal T}=\emptyset$. Let ${\cal I}$
be an initial segment of colex of length $\left|{\cal R}\right|$
in $\left[n\right]^{(i)}$, and set ${\cal J}={\cal S}\setminus{\cal Q}$.
Then $\left|{\cal J}\right|=\left|{\cal R}\right|=\left|{\cal I}\right|$
and ${\cal J}$ is a segment, as ${\cal S}$ and ${\cal Q}$ are initial
segments of colex. Thus Claim 1 implies that $\left|{\cal J}_{1}\right|\leq\left|{\cal I}_{1}\right|$.
On the other hand, Claim 2 implies that $\left|{\cal R}_{1}\right|\geq\left|{\cal I}_{1}\right|$.
Combining these two yields $\left|{\cal R}_{1}\right|\geq\left|{\cal J}_{1}\right|$.
Hence
\[
\left|{\cal S}_{1}\right|+\left|{\cal T}_{1}\right|=\left|{\cal S}_{1}\right|=\left|{\cal J}_{1}\right|+\left|{\cal Q}_{1}\right|\leq\left|{\cal R}_{1}\right|+\left|Q_{1}\right|
\]
as required. \\

\textbf{Case 2. $\left|{\cal Q}\right|+\left|{\cal R}\right|>\left|C_{s}\right|$}

By definition of $B$ it follows that ${\cal S}=\left[n\right]^{(i)}$.
Note that since $\left|{\cal S}\right|\geq\left|{\cal Q}\right|$,
it follows that $\left|{\cal R}\right|\geq\left|{\cal T}\right|$.
Hence ${\cal I}={\cal R}\setminus{\cal T}\subseteq\left[n\right]^{(i+1)}$
is a segment and it satisfies ${\cal R}={\cal I}\cup{\cal T}$. Let
${\cal I}_{*}\subseteq\left[n\right]^{(i+1)}$ be an initial segment
of colex chosen such that ${\cal K}=\left[n\right]^{(i+1)}\setminus{\cal I}_{*}$
is a segment of size $\left|{\cal I}\right|$. Define ${\cal J}=\left[n\right]^{(i)}\setminus{\cal Q}={\cal S}\setminus{\cal Q}$.
Hence ${\cal J}$ is a segment of size $\left|{\cal S}\right|-\left|{\cal Q}\right|=\left|{\cal R}_{1}\right|-\left|{\cal T}_{1}\right|=\left|{\cal I}\right|$.

Claim 3 implies that $\left|{\cal I}_{1}\right|\geq\left|{\cal K}_{1}\right|$
and Claim 4 implies that $\left|{\cal K}_{1}\right|\geq\left|{\cal J}_{1}\right|$.
Thus combining these results yields that $\left|{\cal I}_{1}\right|\geq\left|{\cal J}_{1}\right|$.
Using the definitions of ${\cal I}$ and ${\cal J}$ it follows that
\[
\left|{\cal S}_{1}\right|+\left|{\cal T}_{1}\right|=\left|{\cal J}_{1}\right|+\left|{\cal Q}_{1}\right|+\left|{\cal T}_{1}\right|\leq\left|{\cal I}_{1}\right|+\left|{\cal Q}_{1}\right|+\left|{\cal T}_{1}\right|=\left|{\cal R}_{1}\right|+\left|{\cal Q}_{1}\right|
\]
as required. This completes the proof of Lemma 8. $\square$\\

\textbf{Proof of Lemma 9.}

We will first start by proving Claim 1, and then we will prove that
other claims can be deduced from Claim 1. \\

\textbf{Proof of Claim 1. }Since ${\cal A}$ is a segment, there exists
initial segments ${\cal I}_{A}$ and ${\cal J}_{A}$ of colex with
${\cal A}={\cal I}_{A}\setminus{\cal J}_{A}$, and denote their associated
sets of sequences by $I_{A}$ and $J_{A}$. Let $C$ be obtained from
$J_{A}$ by reversing all the sequences and by adding $2n$ $1$'s
at the start of each reversed sequence. Let $D$ be obtained from
$I$ by adding $2n$ 1's at the end of each sequence in $I$, where
$I$ is the set of sequences associated to ${\cal I}$ . Set $B=C\cup D$. 

Due to the additional 1's at the start of the elements of $C$ and
at the end of the elements in $D$, it follows that $\delta C$ and
$\delta D$ are disjoint sets. Also note that reversing all the sequences
and adding $1$'s to every sequence do not change the size of the
shadow. Hence $\left|\delta B\right|=\left|\delta C\right|+\left|\delta D\right|=\left|\delta I\right|+\left|\delta J_{A}\right|$.
On the other hand, since ${\cal I}$ and ${\cal J}_{A}$ are initial
segments of colex, Lemma 4 implies that $\left|\delta I\right|=\left|{\cal I}_{1}\right|$
and $\left|\delta J_{A}\right|=\left|\left({\cal J}_{A}\right)_{1}\right|$.
Thus 
\begin{equation}
\left|\delta B\right|=\left|{\cal I}_{1}\right|+\left|\left({\cal J}_{A}\right)_{1}\right|.
\end{equation}

Since ${\cal I}_{A}$ is an initial segment of colex, Lemma 4 implies
that $\left|\delta I_{A}\right|=\left|\left({\cal I}_{A}\right)_{1}\right|$.
But ${\cal I}_{A}$ is a disjoint union of ${\cal J}_{A}$ and ${\cal A}$
so 
\begin{equation}
\left|\delta I_{A}\right|=\left|\left({\cal I}_{A}\right)_{1}\right|=\left|\left({\cal J}_{A}\right)_{1}\right|+\left|{\cal A}_{1}\right|.
\end{equation}
Since ${\cal I}_{A}$ is an initial segment of colex, the corresponding
set of sequences $I_{A}$ has minimal shadow inside a connected component.
Since $\left|{\cal B}\right|=\left|{\cal I}_{A}\right|$, it follows
that 
\begin{equation}
\left|\delta B\right|\geq\left|\delta I_{A}\right|.
\end{equation}
Thus combining (4), (5) and (6) yields 
\begin{equation}
\left|{\cal I}_{1}\right|\geq\left|{\cal A}_{1}\right|
\end{equation}
as required. $\square$\\

\textbf{Claim 1 $\Rightarrow$ Claim 3. }Let ${\cal A}$ and ${\cal I}$
be as in Claim 3. Define $\overline{{\cal A}}=\left\{ A^{c}\,:\,A\in{\cal A}\right\} $
and define $\overline{{\cal I}}$ similarly. Note that $\left|\overline{{\cal A}}\right|=\left|{\cal A}\right|$
and $\overline{{\cal A}}\subseteq\left[n\right]^{(n-r)}$. It is easy
to check that if ${\cal B}\subseteq\left[n\right]^{(r)}$ is an initial
segment of colex, then so is $\overline{\left(\left[n\right]^{(r)}\setminus{\cal B}\right)}$
. Thus $\overline{{\cal I}}$ is an initial segment of colex. 

Since ${\cal A}$ is a segment, there exists initial segments ${\cal K}$
and ${\cal L}$ such that ${\cal A}={\cal K}\setminus{\cal L}$. This
can be rewritten as ${\cal A}=\left(\left[n\right]^{\left(r\right)}\setminus{\cal L}\right)\setminus\left(\left[n\right]^{\left(r\right)}\setminus{\cal K}\right)$
and hence $\overline{{\cal A}}=\overline{\left(\left[n\right]^{\left(r\right)}\setminus{\cal L}\right)\setminus\left(\left[n\right]^{\left(r\right)}\setminus{\cal K}\right)}=\overline{\left(\left[n\right]^{\left(r\right)}\setminus{\cal L}\right)}\setminus\overline{\left(\left[n\right]^{\left(r\right)}\setminus{\cal K}\right)}$.
As $\overline{\left(\left[n\right]^{\left(r\right)}\setminus{\cal L}\right)}$
and $\overline{\left(\left[n\right]^{\left(r\right)}\setminus{\cal K}\right)}$
are initial segments of colex, it follows that $\overline{{\cal A}}$
is a segment as well. 

Hence $\overline{{\cal A}}$ and $\overline{{\cal I}}$ satisfies
the conditions of Claim 1, and therefore
\begin{equation}
\left|\left(\overline{{\cal I}}\right)_{1}\right|\geq\left|\left(\overline{{\cal A}}\right)_{1}\right|.
\end{equation}

Note that for any set system ${\cal B}$, we have $\left|{\cal B}\right|=\left|{\cal B}_{1}\right|+\left|\left(\overline{{\cal B}}\right)_{1}\right|$
as for every $A\in{\cal B}$ exactly one of $1\in A$ and $1\in A^{c}$
is satisfied. Thus 
\begin{equation}
\left|{\cal I}\right|=\left|{\cal I}_{1}\right|+\left|\left(\overline{{\cal I}}\right)_{1}\right|
\end{equation}
and 
\begin{equation}
\left|{\cal A}\right|=\left|{\cal A}_{1}\right|+\left|\left(\overline{{\cal A}}\right)_{1}\right|.
\end{equation}
Combining (8), (9) and (10) with $\left|{\cal I}\right|=\left|{\cal A}\right|$
yields that 
\begin{equation}
\left|{\cal A}_{1}\right|\geq\left|{\cal I}_{1}\right|
\end{equation}
which completes the proof of Claim 3.$\square$\\

\textbf{Claim 1 $\Rightarrow$ Claim 2. }Let ${\cal I}$ and ${\cal J}$
be as in Claim 2. For $i+1\leq j\leq n$ let ${\cal S}_{j}=\left\{ A\setminus\{j\}\,:\,A\in{\cal J},\,\max A=j\right\} $.
Thus ${\cal S}_{j}\subseteq\left[j-1\right]^{\left(i\right)}\subseteq\left[n\right]^{\left(i\right)}$
for all $i$. Since ${\cal J}$ is an initial segment of colex, it
follows that ${\cal S}_{j}$ is an initial segment of colex in $\left[j-1\right]^{\left(i\right)}$
for all $j$. ${\cal S}_{j}$ is an initial segment of colex also
in $\left[n\right]^{(i)}$ as initial segments of colex are not affected
by adding new larger elements to the ground set. Note that we can
express ${\cal J}$ as a disjoint union ${\cal J}=\bigcup_{j=i+1}^{n}\left({\cal S}_{j}+\left\{ j\right\} \right)$.
Hence 
\begin{equation}
\left|{\cal J}_{1}\right|=\sum_{j=i+1}^{n}\left|\left({\cal S}_{j}+\left\{ j\right\} \right)_{1}\right|=\sum_{j=i+1}^{n}\left|\left({\cal S}_{j}\right)_{1}\right|.
\end{equation}

Since each ${\cal S}_{j}$ is an initial segment of colex in $\left[n\right]^{(i)}$
and we have$\sum_{j=i+1}^{n}\left|{\cal S}_{j}\right|=\left|{\cal J}\right|=\left|{\cal I}\right|$,
a repeated application of Claim 1 implies that $\left|{\cal J}_{1}\right|\geq\left|{\cal I}_{1}\right|$.
$\square$\\

\textbf{Claim 2 $\Rightarrow$ Claim 4.}

Let ${\cal I}$, ${\cal J}$, ${\cal I}_{*}$ and ${\cal J}_{*}$
be as in Claim 4. Since ${\cal I}_{*}$ and ${\cal J}_{*}$ are initial
segments of colex, the observation pointed out in the proof of Claim
1 $\Rightarrow$ Claim 3 implies that $\overline{{\cal I}}\subseteq\left[n\right]^{(n-i)}$
and $\overline{{\cal J}}\subseteq\left[n\right]^{(n-i-1)}$ are initial
segments of colex as well. Thus Claim 2 implies that 
\begin{equation}
\left|\left(\overline{{\cal I}}\right)_{1}\right|\geq\left|\left(\overline{{\cal J}}\right)_{1}\right|.
\end{equation}
Combining this with 
\begin{equation}
\left|{\cal I}\right|=\left|{\cal I}_{1}\right|+\left|\left(\overline{{\cal I}}\right)_{1}\right|
\end{equation}
and 
\begin{equation}
\left|{\cal J}\right|=\left|{\cal J}_{1}\right|+\left|\left(\overline{{\cal J}}\right)_{1}\right|
\end{equation}
yields $\left|{\cal J}_{1}\right|\geq\left|{\cal I}_{1}\right|$ as
required. $\square$ 

This completes the proof of Lemma 9. $\square$\\

We are now ready to deduce Theorem 1. For convenience, we will recall
the definition of the order $\leq$ and restate Theorem 1. For distinct
$x$, $y\in\left\{ 0,\dots,k\right\} ^{n}$ we set $x\leq y$ if 
\begin{enumerate}
\item $w_{0}\left(x\right)<w_{0}\left(y\right)$ 
\item $w_{0}\left(x\right)=w_{0}\left(y\right)$, $re\left(x\right)\neq re\left(y\right)$
and $re\left(x\right)\leq_{c}re\left(y\right)$
\item $w_{0}\left(x\right)=w_{0}\left(y\right)$, $re\left(x\right)=re\left(y\right)$
and $R_{0}\left(x\right)\leq_{colex}R_{0}\left(y\right)$
\end{enumerate}
\textbf{Theorem 1. }Let $A\subseteq\left\{ 0,\dots,k\right\} ^{n}$
and let $B$ be an initial segment of $\leq$ with $\left|B\right|=\left|A\right|$.
Then $\left|\delta A\right|\geq\left|\delta B\right|$.\\

\textbf{Proof. }Let $A$ be a subset of $\left\{ 0,\dots,k\right\} ^{n}$
of given size with minimal $\delta A$. Define
\[
v(A)=\sum_{j=0}^{n}j\left|A\cap L_{j}\left(n\right)\right|.
\]

If possible, choose $l\in\left[n\right]$, $s\in\left\{ 1,\dots,k\right\} ^{n-l}$
and $t\in\left\{ 1,\dots,k\right\} ^{n-l-1}$ for which $C_{s,t}\left(A\right)\neq A$.
Then by Lemma 8, $B=C_{s,t}\left(A\right)$ satisfies $\left|\delta A\right|\geq\left|\delta B\right|$
and by minimality of $\delta A$ it follows that $\delta B$ is also
minimal. We also have  $v\left(A\right)>v\left(B\right)$, which follows
from the definition of $C_{s,t}\left(A\right)$ and from the fact
that $C_{s,t}\left(A\right)\neq A$. 

Repeating this process we obtain a set $A_{1}$ of size $\left|A\right|$
with minimal $\delta A_{1}$ for which $C_{s,t}\left(A_{1}\right)=A_{1}$
for all $i$, $s\in\left\{ 1,\dots,k\right\} ^{n-l}$ and $t\in\left\{ 1,\dots,k\right\} ^{n-l-1}$.
This follows from the fact that $v\left(B\right)$ is always a non-negative
integer which strictly decreases on each step. Since $C_{s,t}\left(A_{1}\right)=A_{1}$
for all $l\in\left[n\right]$, $s\in\left\{ 1,\dots,k\right\} ^{n-l}$
and $t\in\left\{ 1,\dots,k\right\} ^{n-l-1}$, it is easy to check
that there exists $i$ such that $L_{j}\left(n\right)\subseteq A_{1}$
for all $j<i$ and $L_{j}\left(n\right)\cap A_{1}=\emptyset$ for
all $j>i$. 

Let $C_{s_{1}},\,\dots,\,C_{s_{t}}$ be the connected components in
$L_{i}\left(n\right)$ with $s_{j}\leq_{c}s_{k}$ for $j\leq k$.
Define 
\[
w\left(B\right)=\sum_{j=1}^{t}j\left|C_{s_{j}}\cap B\right|.
\]
If possible, choose $j<k$ for which $C_{s_{j},s_{k}}\left(A_{1}\right)\neq A_{1}$,
and set $B=C_{s_{j},s_{k}}\left(A_{1}\right)$. Now $\left|\delta A_{1}\right|\geq\left|\delta B\right|$
by Lemma 7 and hence $\delta B$ is also minimal. Also $w\left(A_{1}\right)>w\left(B\right)$
follows directly from the definition of the compression operator and
from the definition of $B$. Repeating this process we obtain a set
$A_{2}$ for which 
\begin{enumerate}
\item $\delta A_{2}$ is minimal
\item There exists $i$ such that $L_{j}\left(n\right)\subseteq A_{2}$
for all $j<i$ and $L_{j}\left(n\right)\cap A_{2}=\emptyset$ for
all $j>i$
\item $C_{s_{j},s_{k}}\left(A_{2}\right)=A_{2}$ for all $j<k$
\end{enumerate}
Note that the process must terminate as $w\left(B\right)$ is always
a non-negative integer which strictly decreases on each step. Since
$C_{s_{j},s_{k}}\left(A_{2}\right)=A_{2}$ for all $j<k$ it follows
that there exists $p$ for which $C_{s_{k}}\subseteq A$ for all $k<p$
and $C_{s_{k}}\cap A=\emptyset$ for all $k>p$. 

Let $D=A_{2}\cap C_{s_{p}}$ and let $A_{3}$ be set obtained from
$A_{2}$ by taking $A_{3}\cap C_{s_{p}}$ to be the set corresponding
to an initial segment of colex of length $\left|D\right|$, and taking
$A_{3}\setminus C_{s_{p}}=A_{2}\setminus C_{s_{p}}$. Then Lemma 3
implies that $\left|\delta A_{2}\right|\geq\left|\delta A_{3}\right|$
so $\delta A_{3}$ is minimal. On the other hand, by the construction
of $A_{3}$ it is clear that it is an initial segment of $\leq$.
Hence an initial segment of $\leq$ minimises $\delta$.$\square$

\section{An extremal result for the generalised shadow}

So far we have considered operator which allows us to delete a coordinate
which equals 0. It is natural to ask what happens if we generalise
this set-up and allow the deletion of any coordinate that is in some
chosen set. 

Define $\delta_{r}$-shadow of $A\subseteq\left\{ 0,\dots,k\right\} ^{n}$
to be the subset of sequences in $\left\{ 0,\dots,k\right\} ^{n-1}$
obtained from any of its vectors by removing exactly one coordinate
that is one of $\left\{ 0,\dots,r\right\} $. Thus $\delta=\delta_{0}$
and $\Delta=\delta_{k}$. Define $v_{r}\left(x\right)=\sum_{i=0}^{r}w_{i}(x)$.
That is, $v_{r}\left(x\right)$ is the number of coordinates of $x$
in the set $\left\{ 0,\dots,r\right\} $. Define $L_{s}\left(n\right)=\left\{ x\in\left\{ 0,\dots,k\right\} ^{n}\,:\,v_{r}\left(x\right)=s\right\} $
and $L_{\leq s}\left(n\right)=\bigcup_{i=0}^{s}L_{i}\left(n\right)$.
The aim of this section is to prove that the sets $L_{\leq s}\left(n\right)$
are extremal for $\delta_{r}$. This follows directly from the following
Proposition. \\

\textbf{Proposition 10. }Let $A\subseteq\left\{ 0,\dots,k\right\} ^{n}$
and let $A_{s}=A\cap L_{s}\left(n\right)$. Then 
\[
\left|\delta A\right|\geq\frac{1}{n\left(r+1\right)}\sum_{s=0}^{n}s\left|A_{s}\right|.
\]
\\

\textbf{Proof. }Let $X=\left\{ 0,\dots,k\right\} ^{n}$, $Y=\left\{ 0,\dots,k\right\} ^{n-1}$,
let $H$ be defined as in Section 3 and let ${\cal H}$ be a bipartite
multigraph on $X\cup Y$ with edges given as follows. For each $x\in X\cap L_{s}\left(n\right)$
there are exactly $s$ coordinates $x_{i_{1}},\dots,x_{i_{s}}$ which
are elements of $\left\{ 0,\dots,r\right\} $. Define $y_{j}$ to
be the sequence obtained by deleting the coordinate $x_{i_{j}}$.
Then certainly $y_{j}\in\delta x$ and some of the $y_{j}$ may be
equal. Define the edges of ${\cal H}$ to be the edges $xy_{j}$ for
all $1\leq j\leq s$ counting with multiplicities. For example, when
$r=1$ the sequence $x=00121$ is connected by two edges to $0121$,
and by one edge to both $0012$ and $0021$. 

It is easy to verify that for all $y\in Y$, $y$ has degree $n\left(r+1\right)$
as this corresponds to adding any element of $\left\{ 0,\dots,r\right\} $
to any of the $n$ possible places in the sequence $y$. Note that
for all $x\in X$ we have $\Gamma_{{\cal H}}\left(x\right)=\delta x$,
and hence for any $A\subseteq X$ we have $\delta A=\Gamma_{{\cal H}}\left(A\right)$.
By the definition of ${\cal H}$ we have $d\left(x\right)=s$ for
all $x\in L_{s}\left(n\right)$, and as observed earlier we have $d\left(y\right)=n\left(r+1\right)$
for all $y\in Y$. Since the connected components of ${\cal H}$ are
contained in the sets $L_{s}\left(n\right)\cup L_{s-1}\left(n-1\right)$,
we have $\Gamma_{{\cal H}}\left(A\right)\cap L_{s-1}\left(n-1\right)=\Gamma_{{\cal H}}\left(A\cap L_{s}\left(n\right)\right)$
and therefore 

\begin{equation}
\left|\Gamma_{{\cal H}}A\right|=\sum_{s=0}^{r}\left|\Gamma_{{\cal H}}\left(A_{s}\right)\right|.
\end{equation}

For a set $B\subseteq L_{s}\left(n\right)$ we have 
\[
s\left|B\right|=e\left(B,\Gamma_{{\cal H}}\left(B\right)\right)\leq e\left(\Gamma_{{\cal H}}\left(B\right),X\right)=\left|\Gamma_{{\cal H}}\left(B\right)\right|n\left(r+1\right)
\]
and hence 
\begin{equation}
\left|\Gamma_{{\cal H}}\left(B\right)\right|\geq\frac{s}{n\left(r+1\right)}\left|B\right|.
\end{equation}
Applying (17) to each term of the sum in (16) yields 
\begin{equation}
\left|\delta A\right|=\left|\Gamma_{{\cal H}}A\right|\geq\frac{1}{n\left(r+1\right)}\sum_{s=0}^{r}s\left|A_{s}\right|
\end{equation}
as required. $\square$\\

Now we are ready to conclude that the sets $L_{\leq s}\left(n\right)$
are extremal. \\

\textbf{Corollary 11. }If $A\subseteq\left\{ 0,\dots,k\right\} ^{n}$
and $\left|A\right|=\left|L_{\leq s}\left(n\right)\right|$, then
$\left|\delta A\right|\geq\left|\delta L_{\leq s}\left(n\right)\right|$
with equality if and only if $A=L_{\leq s}\left(n\right)$.\\

\textbf{Proof. }Let $B=L_{\leq s}\left(n\right)$. We will first check
that the equality holds for $B$ in (18). Note that $B_{i}=L_{i}\left(n\right)$
for all $i\leq s$ and $B_{i}=\emptyset$ for all $i>s$. For $i\le s$,
$\left|B_{i}\right|=\left|L_{i}\left(n\right)\right|={n \choose i}\left(r+1\right)^{i}(k-r)^{n-i}$
and $\left|\delta B_{i}\right|=\left|L_{i-1}\left(n-1\right)\right|={n-1 \choose i-1}\left(r+1\right)^{i-1}\left(k-r\right)^{n-i}$.
Therefore $\left|\delta B_{i}\right|=\frac{i}{n\left(r+1\right)}\left|B_{i}\right|$
holds for all $i\leq s$, and in fact also for $i>s$ as in this case
both sides are 0. Hence the equality holds in (17) for all $i$, and
thus the equality holds in (18) as well. 

Given a set $A$ of fixed size with $\left|A_{i}\right|\leq\left|L_{i}\left(n\right)\right|$
for all $i$, it is easy to see that $\frac{1}{n\left(r+1\right)}\sum_{t=0}^{r}t\left|A_{t}\right|$
is minimised if and only if $A=L_{\leq n}\cup B$ for suitably chosen
$n$ and for any $B\subseteq L_{n+1}$ of suitable size. Hence given
$A$ with $\left|A\right|=\left|L_{\leq s}\left(n\right)\right|$,
the quantity $\frac{1}{n\left(r+1\right)}\sum_{t=0}^{r}t\left|A_{t}\right|$
attains its minimum value uniquely when $A=L_{\le s}\left(n\right)$. 

Thus 
\begin{equation}
\left|\delta A\right|\geq\frac{1}{n\left(r+1\right)}\sum_{t=0}^{r}t\left|A_{t}\right|\ge\frac{1}{n\left(r+1\right)}\sum_{t=0}^{s}t\left|L_{t}\left(n\right)\right|=\left|\delta L_{\leq s}\left(n\right)\right|
\end{equation}
and the second inequality holds if and only if $A=L_{\leq s}\left(n\right)$,
as required. $\square$

\end{document}